\documentstyle[amssymb,12pt]{article}

\newtheorem{th}{Theorem}[section]
\newtheorem{lem}[th]{Lemma}

\newtheorem{cor}[th]{Corollary}
\newtheorem{defn}[th]{Definition}
\newenvironment{defn-new}{\begin{defn} \em}{\end{defn}}
\newtheorem{rem}[th]{Remark}
\newenvironment{rem-new}{\begin{rem} \em}{\end{rem}}
\newtheorem{ex}[th]{Example}
\newenvironment{ex-new}{\begin{ex} \em}{\end{ex}}

\makeatletter \@addtoreset{equation}{section} \makeatother
\begin{document}

\bigskip

\begin{center}
{\LARGE {\bf \boldmath $C$-totally real warped product submanifolds}}\bigskip

{\sc Mukut Mani Tripathi}\bigskip \bigskip
\end{center}

\begin{quote}
\noindent {\bf Abstract.} We obtain a basic inequality involving the
Laplacian of the warping function and the squared mean curvature of any
warped product isometrically immersed in a Riemannian manifold (cf. Theorem~%
\ref{th-laplacian-f-H}).
Applying this general theory, we obtain basic inequalities involving the
Laplacian of the warping function and the squared mean curvature of $C$%
-totally real warped product submanifolds of $\left( \kappa ,\mu \right) $%
-space forms, Sasakian space forms and non-Sasakian $\left( \kappa ,\mu
\right) $-manifolds. Then we obtain obstructions to the existence of minimal
isometric immersions of $C$-totally real warped product submanifolds in $%
\left( \kappa ,\mu \right) $-space forms, non-Sasakian $\left( \kappa ,\mu
\right) $-manifolds and Sasakian space forms. In the last, we obtain an
example of a warped product $C$-totally real submanifold of a non-Sasakian $%
\left( \kappa ,\mu \right) $-manifold, which satisfies the equality case of
the basic inequality.


\end{quote}

\section{Introduction\label{sect-prel}}

The concept of warped product was first introduced by Bishop and O'Neill
\cite{Bishop-O'Neill-69} to construct the examples of Riemannian manifolds
with negative curvature. Let $M_{1}$ and $M_{2}$ be two Riemannian manifolds
with Riemannian metrics $g_{1}$ and $g_{2}$ respectively. Let $f>0$ be a
differentiable function on $M_{1}$. Consider the product manifold $M_{1}
\times M_{2}$ with its canonical projections $\pi _{1}:M_{1}\times
M_{2}\rightarrow M_{1}$ and $\pi _{2}:M_{1}\times M_{2}\rightarrow M_{2}$.
The warped product $M=M_{1}\times _{f}M_{2}$ is the product manifold $M_{1}
\times M_{2}$ equipped with the Riemannian structure such that
\begin{equation}
||X||^{2}=||\pi _{1\ast }X||^{2}+\left( f^{2}\circ \pi _{1}\right) \left(
p\right) ||\pi _{2\ast }X)||^{2}  \label{eq-defn-warped-prod-metric-1}
\end{equation}
for each tangent vector $X\in T_{p}M$, $p\in M$, or equivalently, the
Riemannian metric $g$ of $M$ is given by
\begin{equation}
g=g_{1}+f^{2}g_{2}  \label{eq-defn-warped-prod-metric-2}
\end{equation}
with the usual meaning. The function $f$ is called the warping function of
the warped product. Especially, if the warping function is constant, then
the manifold $M$ is said to be {\em trivial}.

\medskip For a warped product $M_{1}\times _{f}M_{2}$ we denote by ${\cal D}%
_{1}$ and ${\cal D}_{2}$ the distributions given by the vectors tangent to
leaves and fibres, respectively. Thus, ${\cal D}_{1}$ is obtained from
tangent vectors to $M_{1}$ via the horizontal lift and ${\cal D}_{2}$ is
obtained by the tangent vectors of $M_{2}$ via the vertical lift. It is well
known that
\begin{equation}
\nabla _{X}Y=\nabla _{Y}X=\frac{1}{f}(Xf)Y,\qquad X\in {\cal D}_{1},\;\;Y\in
{\cal D}_{2}.  \label{eq-del-XY-warp}
\end{equation}
From (\ref{eq-del-XY-warp}) for unit vector fields $X,Z$ tangent to $%
M_{1},M_{2}$, respectively, the sectional curvature $K(X\wedge Z)$ of the
plane section spanned by $X$ and $Z$ becomes
\begin{equation}
K(X\wedge Z)=g(\nabla _{Z}\nabla _{X}X-\nabla _{X}\nabla _{Z}X,Z)=\frac{1}{f}%
\left\{ (\nabla _{X}X)f-X^{2}f\right\} .  \label{eq-K(X,Z)}
\end{equation}
Consequently, we obtain
\begin{equation}
\frac{\Delta f}{f}=\sum_{j=1}^{n_{1}}K\left( e_{j}\wedge e_{s}\right)
,\qquad s\in \{n_{1}+1,\ldots ,n\}.  \label{eq-del-f-upon-f}
\end{equation}

\medskip The notion of warped products plays some important roles in
differential geometry as well as in physics. For instance, the best
relativistic model of the Schwarzschild space-time that describes the outer
space around a massive star or a black hole is given as a warped product
\cite[pp. 364-367]{Neill-book}. For more details we refer to \cite
{Bishop-O'Neill-69} and \cite{Neill-book}. For a survey on warped product as
Riemannian submanifolds we refer to \cite{Chen-02-Soochow}.

\medskip In \cite[Theorem 1.4]{Chen-02-Edinburgh}, B.-Y. Chen established a
sharp relationship between the warping function $f$ of a warped product $%
M_{1}\times _{f}M_{2}$ isometrically immersed in a real space form $R^{m}(c)$
and the squared mean curvature $\Vert H\Vert ^{2}$ given by
\begin{equation}
\frac{\Delta f}{f}\leq \frac{\left( n_{1}+n_{2}\right) ^{2}}{4n_{2}}\Vert
H\Vert ^{2}+n_{1}c,  \label{eq-f-Chen}
\end{equation}
where $n_{i}=\dim \left( M_{i}\right) $, $i=1,2$ and $\bigtriangleup $ is
the Laplacian operator of $M_{1}$. The equality case of (\ref{eq-f-Chen})
holds identically if and only if $x$ is a mixed totally geodesic immersion
and $n_{1}H_{1}=n_{2}H_{2}$, where $H_{i}$, $i=1,2,$ are the partial mean
curvature vectors. As a refinement, in \cite[Theorem 5]{Chen-05-JGeom}, he
proved that the equality case of (\ref{eq-f-Chen}) is true if and only if
one of the following two cases occurs: (1) The warping function $f$ is an
eigenfunction of the Laplacian operator $\bigtriangleup $ with eigenvalue $%
n_{1}c$ and $x$ is a minimal immersion; (2) $\bigtriangleup f\neq \left(
n_{1}c\right) f$ and locally $x$ is a non-minimal warped product immersion $%
\left( x_{1},x_{2}\right) :M_{1}\times _{f}M_{2}\rightarrow N_{1}\times
_{\rho }N_{2}$ of $M_{1}\times _{f}M_{2}$ into some warped product
representation $N_{1}\times _{\rho }N_{2}$ of the real space form $\tilde{M}%
(c)$ such that $x_{2}:M_{2}\rightarrow N_{2}$ is a minimal immersion and the
mean curvature vector of $x_{1}:M_{1}\rightarrow N_{1}$ is given by $-\left(
n_{2}/n_{1}\right) D\ln \rho $.

\medskip The inequality (\ref{eq-f-Chen}) was noticed by several authors and
they established same kind of inequalities for different kind of
submanifolds in ambient manifolds possessing different kind of structures
(for example, see \cite{Chen-02-PJA}, \cite{Chen-03-PJA}, \cite
{Chen-04-Rocky}, \cite{Kim-Yoon-04}, \cite{Mat-Mih-02-SUT}, \cite{MihA-04},
\cite{MihA-05}, \cite{Yoon-04}, \cite{YC-04}, \cite{YCH-05}). Now, it is a
natural motivation to find a basic inequality involving the warping function
and the squared mean curvature of any warped product isometrically immersed
in any Riemannian manifold without assuming any restriction on the Riemann
curvature tensor of the ambient manifold. Using technique of B.-Y. Chen and
the concept of the scalar curvature of $k$-plane sections, in section~\ref
{sect-basic-ineq-warp-prod-submfd} we achieve this goal by obtaining a basic
inequality involving the Laplacian of the warping function $f$ and the
squared mean curvature of a warped product $M_{1}\times _{f}M_{2}$
isometrically immersed in a Riemannian manifold (see Theorem~\ref
{th-laplacian-f-H}). Section~\ref{sect-kmu} contains some necessary
background of contact geometry including the concepts of Sasakian manifolds,
$\left( \kappa ,\mu \right) $-manifolds, $\left( \kappa ,\mu \right) $-space
forms and non-Sasakian $\left( \kappa ,\mu \right) $-manifolds. In section~%
\ref{sect-C-tot-real-wp-submfd}, we apply the general theory given by
Theorem~\ref{th-laplacian-f-H} to obtain corresponding results for $C$%
-totally real warped product submanifolds of $\left( \kappa ,\mu \right) $%
-space forms, Sasakian space forms and non-Sasakian $\left( \kappa ,\mu
\right) $-manifolds. Then we obtain obstructions to the existence of minimal
isometric immersions of $C$-totally real warped product submanifolds in
these spaces. In the last, we also obtain an example of a warped product $C$%
-totally real submanifold of a non-Sasakian $\left( \kappa ,\mu \right) $%
-manifold which satisfies the equality case of the corresponding basic
inequality.

\section{A basic inequality for warped product submanifolds\label%
{sect-basic-ineq-warp-prod-submfd}}

Let $M$ be an $n$-dimensional Riemannian manifold equipped with a Riemannian
metric $g$. The inner product of the metric $g$ is denoted by $\left\langle
\,,\right\rangle$. Let $\left\{ e_{1},\ldots ,e_{n}\right\} $ be any
orthonormal basis for $T_{p}M$. The {\em scalar curvature} $\tau \left(
p\right) $ of $M$ at $p$ is defined by
\begin{equation}
\tau \left( p\right) =\sum_{1\leq i<j\leq n}K\left( e_{i}\wedge e_{j}\right),
\label{eq-tau(p)}
\end{equation}
where $K\left( e_{i}\wedge e_{j}\right) $ is the {\em sectional curvature}
of the plane section spanned by $e_{i}$ and $e_{j}$ at $p\in M$. Let $\Pi
_{k}$ be a $k$-plane section of $T_{p}M$ and $\{e_{1},\ldots ,e_{k}\}$ any
orthonormal basis of $\Pi _{k}$. The scalar curvature $\tau \left( \Pi
_{k}\right) $ of $\Pi _{k}$ is given by \cite{Chen-99-GMJ}
\begin{equation}
\tau \left( \Pi _{k}\right) =\sum_{1\leq i<j\leq k}K\left( e_{i}\wedge
e_{j}\right) .  \label{eq-tau(Pi-k)}
\end{equation}
The scalar curvature $\tau (p)$ of $M$ at $p$ is identical with the scalar
curvature of the tangent space $T_{p}M$ of $M$ at $p$, that is, $\tau \left(
p\right) =\tau \left( T_{p}M\right) $. Geometrically, $\tau (\Pi _{k})$ is
the scalar curvature of the image ${\rm {exp}}_{p}(\Pi _{k})$ of $\Pi _{k}$
at $p$ under the exponential map at $p$. If $\Pi _{2}$ is a $2$-plane
section, $\tau (\Pi _{2})$ is simply the sectional curvature $K(\Pi_{2})$ of
$\Pi_{2}$. \medskip

Let $(M,g)$ be a submanifold of a Riemannian manifold $\tilde{M}$ equipped
with a Riemannian metric $\tilde{g}$. Covariant derivatives and curvatures
with respect to $\left( M,g\right) $ will be written in the usual manner,
while those with respect to the ambient manifold $(\tilde{M}, \tilde{g}) $
will be written with ``tildes'' over them. We use the inner product notation
$\left\langle\; ,\right\rangle $ for the metric $\tilde{g}$ of $\tilde{M}$
as well as for the induced metrics on the submanifold $M$ and its normal
bundle.

\medskip The Gauss and Weingarten formulas are given respectively by
\[
\tilde{\nabla }_{X}Y=\nabla _{X}Y+\sigma \left( X,Y\right) \qquad {\rm and}%
\qquad \tilde{\nabla }_{X}N=-A_{N}X+\nabla _{X}^{\perp }N
\]
for all $X,Y\in TM$ and $N\in T^{\perp }M$, where $\tilde{\nabla }$, $\nabla$
and $\nabla^{\perp}$ are respectively the Riemannian, induced Riemannian and
induced normal connections in $\tilde{M}$, $M$ and the normal bundle $%
T^{\perp }M$ of $M$ respectively, and $\sigma $ is the second fundamental
form related to the shape operator $A$ by $\left\langle \sigma \left(
X,Y\right) , N\right\rangle = \left\langle A_{N}X,Y\right\rangle $. The
equation of Gauss is given by
\begin{eqnarray}
R(X,Y,Z,W) &=&\tilde{R}(X,Y,Z,W)+\ \left\langle \sigma (X,W),\sigma
(Y,Z)\right\rangle  \nonumber \\
&&-\ \left\langle \sigma (X,Z),\sigma (Y,W)\right\rangle  \label{eq-Gauss}
\end{eqnarray}
for all $X,Y,Z,W\in TM$, where $\tilde{R}$ and $R$ are the curvature tensors
of $\tilde{M}$ and $M$ respectively.

\medskip For any orthonormal basis $\left\{ e_{1},\ldots ,e_{n}\right\}$ of
the tangent space $T_{p}M$, the mean curvature vector $H\left( p\right) $ is
given by
\begin{equation}
H\left( p\right) =\frac{1}{n}\sum_{i=1}^{n}\sigma \left( e_{i},e_{i}\right) ,
\label{eq-mean-curv}
\end{equation}
where $n={\rm dim}(M)$. The submanifold $M$ is {\em totally geodesic} in $%
\tilde{M}$ if $\sigma =0$, and {\em minimal} if $H=0$. If $\sigma \left(
X,Y\right) $ $=$ $g\left( X,Y\right) H$ for all $X,Y\in TM$, then $M$ is
{\em totally umbilical}.

\medskip Now, let $\left\{ e_{1},\ldots ,e_{n}\right\} $ be an orthonormal
basis of the tangent space $T_{p}M$ and $e_{r}$ belongs to an orthonormal
basis $\left\{ e_{n+1},\ldots ,e_{m}\right\} $ of the normal space $%
T_{p}^{\perp }M$. We put
\[
\sigma _{ij}^{r}=\left\langle \sigma \left( e_{i},e_{j}\right)
,e_{r}\right\rangle \quad {\rm and}\quad \left\Vert \sigma \right\Vert
^{2}=\sum_{i,j=1}^{n}\left\langle \sigma \left( e_{i},e_{j}\right) ,\sigma
\left( e_{i},e_{j}\right) \right\rangle .
\]
Let $K\left( e_{i}\wedge e_{j}\right) $ and $\tilde{K}\left( e_{i}\wedge
e_{j}\right) $ denote the sectional curvature of the plane section spanned
by $e_{i}$ and $e_{j}$ at $p$ in the submanifold $M$ and in the ambient
manifold $\tilde{M}$ respectively. In view of the equation (\ref{eq-Gauss})
of Gauss, we have
\begin{equation}
K\left( e_{i}\wedge e_{j}\right) =\tilde{K}\left( e_{i}\wedge e_{j}\right)
+\sum_{r=n+1}^{m}\left( \sigma _{ii}^{r}\sigma _{jj}^{r}-(\sigma
_{ij}^{r})^{2}\right) .  \label{eq-Kij}
\end{equation}
From (\ref{eq-Kij}) it follows that
\begin{equation}
2\tau \left( p\right) =2\tilde{\tau}\left( T_{p}M\right) +n^{2}\left\Vert
H\right\Vert ^{2}-\left\Vert \sigma \right\Vert ^{2},  \label{eq-tau-H-sigma}
\end{equation}
where
\[
\tilde{\tau}\left( T_{p}M\right) =\sum_{1\leq i<j\leq n}\tilde{K}\left(
e_{i}\wedge e_{j}\right)
\]
denote the scalar curvature of the $n$-plane section $T_{p}M$ in the ambient
manifold $\tilde{M}$.

\medskip Now, let $x:M_{1}\times _{f}M_{2}\rightarrow \tilde{M}$ be an
isometric immersion of a warped product $M_{1}\times _{f}M_{2}$ into a
Riemannian manifold $\tilde{M}$. We write $n_{i}H_{i}={\rm trace}\,\left(
\sigma _{i}|_{M_{i}}\right) $, where trace$\left( \sigma
_{i}|_{M_{i}}\right) $ is the trace of the second fundamental form $\sigma $
of $x$ restricted to $M_{i}$ and $n_{i}=\dim M_{i}$ $(i=1,2)$. The immersion
$x$ is called {\em mixed totally geodesic} if $\sigma \left( X,Y\right) =0$
for $X\in {\cal D}_{1}$ and $Y\in {\cal D}_{2}$.

\medskip We shall need the following Lemma.

\begin{lem}
\label{lem-Chen} {\bf \cite{Chen-93}} Let $\ell \geq 2$ and $a_{1},\ldots
,a_{\ell }\,,b$ be real numbers such that
\[
\left( \sum_{i=1}^{\ell }a_{i}\right) ^{2}=(\ell -1)\left( \sum_{i=1}^{\ell
}a_{i}^{2}+b\right)
\]
Then $2a_{1}a_{2}\geq b$, with equality holding if and only if
\[
a_{1}+a_{2}=a_{3}=\cdots =a_{\ell }.
\]
\end{lem}

Now, we obtain a basic inequality involving the Laplacian of the warping
function and the squared mean curvature of a warped product submanifold of a
Riemannian manifold.

\begin{th}
\label{th-laplacian-f-H} Let $x$ be an isometric immersion of an $n$%
-dimensional warped product manifold $M_{1}\times _{f}M_{2}$ into an $m$%
-dimensional Riemannian manifold $\tilde{M}$. Then
\begin{equation}
n_{2}\frac{\Delta f}{f}\leq \frac{n^{2}}{4}\Vert H\Vert ^{2}+\tilde{\tau}%
\left( T_{p}M\right) -\tilde{\tau}\left( T_{p}M_{1}\right) -\tilde{\tau}%
\left( T_{p}M_{2}\right) ,  \label{eq-lap-f-H}
\end{equation}
where $n_{i}=\dim M_{i}$, $i=1,2$, and $\Delta $ is the Laplacian operator
of $M_{1}$.

Moreover, the equality case of $(\ref{eq-lap-f-H})$ holds identically if and
only if $x$ is a mixed totally geodesic immersion and $n_{1}H_{1}=n_{2}H_{2}$%
, where $H_{i}$, $i=1,2,$ are the partial mean curvature vectors.
\end{th}

\noindent {\bf Proof.} Let $M=M_{1}\times _{f}M_{2}$ be a warped product
submanifold of a Riemannian manifold $\tilde{M}$. We set
\begin{equation}
2\delta =4\tau \left( p\right) -4\tilde{\tau}\left( T_{p}M\right)
-n^{2}\Vert H\Vert ^{2},  \label{eq-delta-def}
\end{equation}
so that the equation (\ref{eq-tau-H-sigma}) can be written as
\begin{equation}
n^{2}\Vert H\Vert ^{2}=2(\delta +\Vert \sigma \Vert ^{2}).
\label{eq-H-delta-sigma}
\end{equation}
Choose a local orthonormal frame $\{e_{1},\ldots ,e_{n},e_{n+1},\ldots
,e_{m}\}$, such that $e_{1},\ldots ,e_{n_{1}}$ are tangent to $M_{1}$, $%
e_{n_{1}+1},\ldots ,e_{n}$ are tangent to $M_{2}$ and $e_{n+1}$ is parallel
to the mean curvature vector $H$. If we put
\[
a_{1}=\sigma _{11}^{n+1},\qquad a_{2}=\sum_{i=2}^{n_{1}}\sigma
_{ii}^{n+1},\qquad a_{3}=\sum_{t=n_{1}+1}^{n}\sigma _{tt}^{n+1},
\]
then with respect to the above orthonormal frame the equation (\ref
{eq-H-delta-sigma}) becomes
\begin{equation}
\left( \sum_{i=1}^{3}a_{i}\right) ^{2}=2\left(
\sum_{i=1}^{3}a_{i}^{2}+b\right) ,  \label{eq-a-i}
\end{equation}
where
\begin{eqnarray}
b &=&\left\{ \delta +\sum_{1\leq i\neq j\leq n}(\sigma
_{ij}^{n+1})^{2}+\sum_{r=n+2}^{m}\sum_{i,j=1}^{n}(\sigma
_{ij}^{r})^{2}\right.  \nonumber \\
&&\qquad \left. -\ \sum_{2\leq j\neq k\leq n_{1}}\sigma _{jj}^{n+1}\sigma
_{kk}^{n+1}-\sum_{n_{1}+1\leq s\neq t\leq n}\sigma _{ss}^{n+1}\sigma
_{tt}^{n+1}\right\} .  \label{eq-b}
\end{eqnarray}
Applying Lemma~\ref{lem-Chen} to (\ref{eq-a-i}), we get $2a_{1}a_{2}\geq b$,
with equality holding if and only if $a_{1}+a_{2}=a_{3}$. Equivalently, we
get
\begin{eqnarray}
&&\sum_{1\leq j<k\leq n_{1}}\sigma _{jj}^{n+1}\sigma
_{kk}^{n+1}+\sum_{n_{1}+1\leq s<t\leq n}\sigma _{ss}^{n+1}\sigma _{tt}^{n+1}
\nonumber \\
&&\quad \geq \quad \frac{\delta }{2}+\sum_{1\leq \alpha <\beta \leq
n}(\sigma _{\alpha \beta }^{n+1})^{2}+\frac{1}{2}\sum_{r=n+2}^{m}\sum_{%
\alpha ,\beta =1}^{n}(\sigma _{\alpha \beta }^{r})^{2}  \label{eq-a1-a2-b}
\end{eqnarray}
with equality holding if and only if
\begin{equation}
\sum_{i=1}^{n_{1}}\sigma _{ii}^{n+1}=\sum_{s=n_{1}+1}^{n}\sigma _{ss}^{n+1}.
\label{eq-a1-a2-a3}
\end{equation}

From equation (\ref{eq-del-f-upon-f}) we get
\[
n_{2}\frac{\Delta f}{f}=\sum_{j=1}^{n_{1}}\sum_{s=n_{1}+1}^{n}K\left(e_{j}%
\wedge e_{s}\right) =\tau \left( p\right) -\tau \left( T_{p}M_{1}\right)
-\tau \left( T_{p}M_{2}\right) ,
\]
which in view of (\ref{eq-Gauss}) gives
\begin{eqnarray}
n_{2}\frac{\Delta f}{f} &=&\tau \left( p\right) -\tilde{\tau}\left(
T_{p}M_{1}\right) -\sum_{r=n+1}^{m}\sum_{1\leq j<k\leq n_{1}}(\sigma
_{jj}^{r}\sigma _{kk}^{r}-(\sigma _{jk}^{r})^{2})  \nonumber \\
&&-\ \tilde{\tau}\left( T_{p}M_{2}\right) \sum_{r=n+1}^{m}\sum_{n_{1}+1\leq
s<t\leq n}(\sigma _{ss}^{r}\sigma _{tt}^{r}-(\sigma _{st}^{r})^{2})\,.
\label{eq-del-f-upon-f-2}
\end{eqnarray}

In view of (\ref{eq-a1-a2-b}), (\ref{eq-del-f-upon-f-2}) and (\ref
{eq-delta-def}) we get
\begin{eqnarray*}
n_{2}\frac{\Delta f}{f} &\leq &\frac{n^{2}}{4}\Vert H\Vert ^{2}+\tilde{\tau}%
\left( T_{p}M\right) -\tilde{\tau}\left( T_{p}M_{1}\right) -\tilde{\tau}%
\left( T_{p}M_{2}\right) \\
&&-\ \sum_{j=1}^{n_{1}}\sum_{t=n_{1}+1}^{n}(\sigma _{jt}^{n+1})^{2}-\frac{1}{%
2}\sum_{r=n+2}^{m}\sum_{\alpha ,\beta =1}^{n}(\sigma _{\alpha \beta
}^{r})^{2} \\
&&+\ \sum_{r=n+2}^{m}\sum_{1\leq j<k\leq n_{1}}((\sigma
_{jk}^{r})^{2}-\sigma _{jj}^{r}\sigma _{kk}^{r})+\
\sum_{r=n+2}^{m}\sum_{n_{1}+1\leq s<t\leq n}((\sigma _{st}^{r})^{2}-\sigma
_{ss}^{r}\sigma _{tt}^{r})
\end{eqnarray*}
or
\begin{eqnarray}
n_{2}\frac{\Delta f}{f} &\leq &\frac{n^{2}}{4}\Vert H\Vert ^{2}+\tilde{\tau}%
\left( T_{p}M\right) -\tilde{\tau}\left( T_{p}M_{1}\right) -\tilde{\tau}%
\left( T_{p}M_{2}\right)
-\sum_{r=n+1}^{m}\sum_{j=1}^{n_{1}}\sum_{t=n_{1}+1}^{n}(\sigma _{jt}^{r})^{2}
\nonumber \\
&&-\frac{1}{2}\sum_{r=n+2}^{m}\left( \sum_{j=1}^{n_{1}}\sigma
_{jj}^{r}\right) ^{2}-\frac{1}{2}\sum_{r=n+2}^{m}\left(
\sum_{t=n_{1}+1}^{n}\sigma _{tt}^{r}\right) ^{2},  \label{eq-del-f-upon-f-3}
\end{eqnarray}
which implies the inequality (\ref{eq-lap-f-H}). \medskip

The equality sign of (\ref{eq-del-f-upon-f-3}) is true if and only if
\begin{equation}
\sigma _{jt}^{r}=0,\qquad 1\leq j\leq n_{1},\;n_{1}+1\leq t\leq n,\;n+1 \leq
r\leq m,  \label{eq-mixed-totally-geod}
\end{equation}
and
\begin{equation}
\sum_{i=1}^{n_{1}}\sigma _{ii}^{r}=0=\sum_{t=n_{1}+1}^{n}\sigma
_{tt}^{r},\qquad n+2\leq r\leq m.  \label{eq-n1-H1=n2-H2}
\end{equation}
Obviously (\ref{eq-mixed-totally-geod}) is true if and only if the warped
product $M_{1}\times _{f}M_{2}$ is mixed totally geodesic. From the
equations (\ref{eq-a1-a2-a3}) and (\ref{eq-n1-H1=n2-H2}) it follows that $%
n_{1}H_{1}=n_{2}H_{2}$.

\medskip The converse statement is straightforward. $\blacksquare$

\section{\protect\boldmath $\left( \kappa ,\mu \right) $-manifolds \label%
{sect-kmu}}

A $\left( 2m+1\right) $-dimensional differentiable manifold $\tilde{M}$ is
called an almost contact metric manifold if there is an almost contact
metric structure $(\varphi ,\xi ,\eta ,g)$ consisting of a $\left(
1,1\right) $ tensor field $\varphi $, a vector field $\xi $, a $1$-form $%
\eta $ and a compatible Riemannian metric $g$ satisfying
\begin{equation}
\varphi ^{2}=-I+\eta \otimes \xi ,\quad \eta (\xi )=1,\quad \varphi \xi
=0,\quad \eta \circ \varphi =0,  \label{eq-phi-eta-xi}
\end{equation}
\begin{equation}
\left\langle X,Y\right\rangle =\left\langle \varphi X,\varphi Y\right\rangle
+\eta (X)\eta (Y),  \label{eq-metric-1}
\end{equation}
\begin{equation}
\left\langle X,\varphi Y\right\rangle =-\,\left\langle \varphi
X,Y\right\rangle ,\quad \left\langle X,\xi \right\rangle =\eta (X)
\label{eq-metric-2}
\end{equation}
for all $X,Y\in T\tilde{M}$, where $\left\langle \;,\right\rangle $ the
inner product of the metric $g$. An almost contact metric structure becomes
a {\em contact metric structure} if
\[
\langle X,\varphi Y\rangle =d\eta (X,Y),\qquad X,Y\in T\tilde{M}.
\]
A contact metric structure is a {\em Sasakian\ structure }if and only if the
Riemann curvature tensor $\tilde{R}$ satisfies
\begin{equation}
\tilde{R}(X,Y)\xi =\eta (Y)X-\eta (X)Y,\qquad X,Y\in T\tilde{M}.
\label{eq-Sas-2}
\end{equation}
In a contact metric manifold $\tilde{M}$, the $\left( 1,1\right) $-tensor
field $h$ defined by $2h={\frak L}_{\xi }\varphi $, which is the Lie
derivative of $\varphi $ in the characteristic direction $\xi $, is
symmetric and satisfies
\begin{equation}
h\xi =0,\;\;h\varphi +\varphi h=0,\;\;\tilde{\nabla}\xi =-\varphi -\varphi
h,\;\;{\rm trace}(h)={\rm trace}(\varphi h)=0,  \label{eq-cont-h}
\end{equation}
where $\tilde{\nabla}$ is Levi-Civita connection. In \cite{BKP-95}, D. E.
Blair, T. Koufogiorgos and B. J. Papantoniou introduced the class of contact
metric manifolds $\tilde{M}$ with contact metric structures $\left( \varphi
,\xi ,\eta ,g\right) $, which satisfy
\begin{equation}
\tilde{R}(X,Y)\xi =\left( \kappa I+\mu h\right) \left( \eta \left( Y\right)
X-\eta \left( X\right) Y\right) ,\qquad X,Y\in T\tilde{M},  \label{eq-km}
\end{equation}
where $\kappa $, $\mu $ are real constants. A contact metric manifold
belonging to this class is called a $\left( \kappa ,\mu \right) ${\em %
-manifold}. In a $(\kappa ,\mu )$-manifold we have $h^{2}=\left( \kappa
-1\right) \varphi ^{2}$ and $\kappa \leq 1$. For a $(\kappa ,\mu )$%
-manifold, the conditions of being a Sasakian manifold, $\kappa =1$ and $h=0$
are all equivalent. If $\kappa <1$, and the eigenvalues of $h$ are $0$, $%
\lambda $ and $-\lambda $, where $\lambda =\sqrt{1-\kappa }$. The eigenspace
relative to the eigenvalue $0$ is $\left\{ \xi \right\} $. Moreover, for $%
\kappa \neq 1$, the subbundle ${\cal D}={\rm ker}(\eta )$ can be decomposed
in the eigenspace distributions ${\cal D}_{+}$ and ${\cal D}_{-}$ relative
to the eigenvalues $\lambda $ and $-\lambda $, respectively. These
distributions are orthogonal to each other and have dimension $n$. There are
many motivations for studying $\left( \kappa ,\mu \right) $-manifolds: the
first is that, in the non-Sasakian case the condition (\ref{eq-km})
determines the curvature completely; moreover, while the values of $\kappa $
and $\mu $ change, the form of (\ref{eq-km}) is invariant under $D$%
-homothetic deformations (\cite{BKP-95}); finally, there is a complete
classification of these manifolds, given in \cite{Boeckx-00} by E. Boeckx,
who proved also that any non-Sasakian $(\kappa ,\mu )$-manifold is locally
homogeneous and strongly locally $\varphi $-symmetric (\cite{Boeckx-99},
\cite{Boeckx-06}). There are also non-trivial examples of $(\kappa ,\mu )$%
-manifolds, the most important being the unit tangent sphere bundle $T_{1}%
\tilde{M}$ of a Riemannian manifold $\tilde{M}$ of constant sectional
curvature with the usual contact metric structure. In particular if $\tilde{M%
}$ has constant curvature $c$, then $\kappa =c(2-c)$ and $\mu =-2c$.
Characteristic examples of non-Sasakian $(\kappa ,\mu )$-manifolds are the
tangent sphere bundles of Riemannian manifolds of constant sectional
curvature not equal to one and certain Lie groups \cite{Boeckx-00}. For more
details we refer to \cite{Blair-book}, \cite{BKP-95} and \cite{Koufog-97}.

\medskip Like complex space forms in Hermitian geometry, in contact geometry
we have the notion of manifolds with constant $\varphi $-sectional
curvature. In an almost contact metric manifold, for a unit vector $X$
orthogonal to $\xi $, the sectional curvature $\tilde{K}(X\wedge \varphi X)$
is called a $\varphi $-sectional curvature. In \cite{Koufog-97}, T.
Koufogiorgos showed that in a $(\kappa ,\mu )$-manifold $(\tilde{M},\varphi
,\xi ,\eta ,g)$ of dimension $>3$ if the $\varphi $-sectional curvature at a
point $p$ is independent of the $\varphi $-section at $p$, then it is
constant. If the $\left( \kappa ,\mu \right) $-manifold $\tilde{M}$ has
constant $\varphi $-sectional curvature $c$ then it is called a $\left(
\kappa ,\mu \right) ${\em -space form} and is denoted by $\tilde{M}(c)$. The
Riemann curvature tensor $\tilde{R}$ of $\tilde{M}(c)$ is given explicitly
in its $\left( 0,4\right) $-form by \cite{Koufog-97}
\begin{eqnarray}
&&\tilde{R}\left( X,Y,Z,W\right) \ =\ \frac{c+3}{4}\left\{ \left\langle
Y,Z\right\rangle \left\langle X,W\right\rangle -\left\langle
X,Z\right\rangle \left\langle Y,W\right\rangle \right\}   \nonumber \\
&&+\ \frac{c-1}{4}\{2\left\langle X,\varphi Y\right\rangle \left\langle
\varphi Z,W\right\rangle +\left\langle X,\varphi Z\right\rangle \left\langle
\varphi Y,W\right\rangle -\left\langle Y,\varphi Z\right\rangle \left\langle
\varphi X,W\right\rangle \}  \nonumber \\
&&+\ \frac{c+3-4\kappa }{4}\left\{ \eta \left( X\right) \eta \left( Z\right)
\left\langle Y,W\right\rangle -\eta \left( Y\right) \eta \left( Z\right)
\left\langle X,W\right\rangle \right.   \nonumber \\
&&\qquad \qquad \qquad \left. +\ \left\langle X,Z\right\rangle \eta \left(
Y\right) \eta \left( W\right) -\left\langle Y,Z\right\rangle \eta \left(
X\right) \eta \left( W\right) \right\}   \nonumber \\
&&+\ \frac{1}{2}\left\{ \left\langle hY,Z\right\rangle \left\langle
hX,W\right\rangle -\left\langle hX,Z\right\rangle \left\langle
hY,W\right\rangle \right.   \nonumber \\
&&\qquad \left. +\ \left\langle \varphi hX,Z\right\rangle \left\langle
\varphi hY,W\right\rangle -\left\langle \varphi hY,Z\right\rangle
\left\langle \varphi hX,W\right\rangle \right\}   \nonumber \\
&&+\ \left\langle \varphi Y,\varphi Z\right\rangle \left\langle
hX,W\right\rangle -\left\langle \varphi X,\varphi Z\right\rangle
\left\langle hY,W\right\rangle   \nonumber \\
&&+\ \left\langle hX,Z\right\rangle \left\langle \varphi
^{2}Y,W\right\rangle -\left\langle hY,Z\right\rangle \left\langle \varphi
^{2}X,W\right\rangle   \nonumber \\
&&+\ \mu \left\{ \eta \left( Y\right) \eta \left( Z\right) \left\langle
hX,W\right\rangle -\eta \left( X\right) \eta \left( Z\right) \left\langle
hY,W\right\rangle \right.   \nonumber \\
&&\qquad \left. +\ \left\langle hY,Z\right\rangle \eta \left( X\right) \eta
\left( W\right) -\left\langle hX,Z\right\rangle \eta \left( Y\right) \eta
\left( W\right) \right\}   \label{eq-Riem-curv}
\end{eqnarray}
for all $X,Y,Z,W\in T\tilde{M}$. In particular, if $\kappa =1$ then a $%
\left( \kappa ,\mu \right) $-space form $\tilde{M}(c)$ reduces to a Sasakian
space form $\tilde{M}(c)$ and (\ref{eq-Riem-curv}) reduces to
\begin{eqnarray}
&&\tilde{R}\left( X,Y,Z,W\right) \ =\ \frac{c+3}{4}\left\{ \left\langle
Y,Z\right\rangle \left\langle X,W\right\rangle -\left\langle
X,Z\right\rangle \left\langle Y,W\right\rangle \right\}   \nonumber \\
&&+\,\frac{c-1}{4}\left\{ 2\left\langle X,\varphi Y\right\rangle
\left\langle \varphi Z,W\right\rangle +\left\langle X,\varphi Z\right\rangle
\left\langle \varphi Y,W\right\rangle -\left\langle Y,\varphi Z\right\rangle
\left\langle \varphi X,W\right\rangle \right.   \nonumber \\
&&\left. \qquad +\ \eta \left( X\right) \eta \left( Z\right) \left\langle
Y,W\right\rangle -\eta \left( Y\right) \eta \left( Z\right) \left\langle
X,W\right\rangle \right.   \nonumber \\
&&\left. \qquad +\ \left\langle X,Z\right\rangle \eta \left( Y\right) \eta
\left( W\right) -\left\langle Y,Z\right\rangle \eta \left( X\right) \eta
\left( W\right) \right\} .
\end{eqnarray}
For a non-Sasakian $\left( \kappa ,\mu \right) $-manifold $\tilde{M}$, its
Riemann curvature tensor $\tilde{R}$ is given explicitly in its $\left(
0,4\right) $-form by (\cite{Boeckx-99},\cite{Boeckx-00})
\begin{eqnarray}
&&\tilde{R}\left( X,Y,Z,W\right) =\left( 1-\frac{\mu }{2}\right) \left\{
\left\langle Y,Z\right\rangle \left\langle X,W\right\rangle -\left\langle
X,Z\right\rangle \left\langle Y,W\right\rangle \right\}   \nonumber \\
&&\qquad -\,\frac{\mu }{2}\left\{ 2\left\langle X,\varphi Y\right\rangle
\left\langle \varphi Z,W\right\rangle +\left\langle X,\varphi Z\right\rangle
\left\langle \varphi Y,W\right\rangle -\left\langle Y,\varphi Z\right\rangle
\left\langle \varphi X,W\right\rangle \right\}   \nonumber \\
&&\qquad +\,\left\langle Y,Z\right\rangle \left\langle hX,W\right\rangle
-\left\langle X,Z\right\rangle \left\langle hY,W\right\rangle
-\,\left\langle Y,W\right\rangle \left\langle hX,Z\right\rangle
+\left\langle X,W\right\rangle \left\langle hY,Z\right\rangle   \nonumber \\
&&\qquad +\,\frac{1-\left( \mu /2\right) }{1-\kappa }\left\{ \left\langle
hY,Z\right\rangle \left\langle hX,W\right\rangle -\left\langle
hX,Z\right\rangle \left\langle hY,W\right\rangle \right\}   \nonumber \\
&&\qquad +\,\frac{\kappa -\left( \mu /2\right) }{1-\kappa }\left\{
\left\langle \varphi hY,Z\right\rangle \left\langle \varphi
hX,W\right\rangle -\left\langle \varphi hX,Z\right\rangle \left\langle
\varphi hY,W\right\rangle \right\}   \nonumber \\
&&\qquad +\,\eta \left( X\right) \eta \left( W\right) \left\{ \left( \kappa
-1+\left( \mu /2\right) \right) \left\langle Y,Z\right\rangle +\left( \mu
-1\right) \left\langle hY,Z\right\rangle \right\}   \nonumber \\
&&\qquad -\,\eta \left( X\right) \eta \left( Z\right) \left\{ \left( \kappa
-1+\left( \mu /2\right) \right) \left\langle Y,W\right\rangle +\left( \mu
-1\right) \left\langle hY,W\right\rangle \right\}   \nonumber \\
&&\qquad +\,\eta \left( Y\right) \eta \left( Z\right) \left\{ \left( \kappa
-1+\left( \mu /2\right) \right) \left\langle X,W\right\rangle +\left( \mu
-1\right) \left\langle hX,W\right\rangle \right\}   \nonumber \\
&&\qquad -\,\eta \left( Y\right) \eta \left( W\right) \left\{ \left( \kappa
-1+\left( \mu /2\right) \right) \left\langle X,Z\right\rangle +\left( \mu
-1\right) \left\langle hX,Z\right\rangle \right\}   \label{eq-kmu-curvature}
\end{eqnarray}
for all vector fields $X,Y,Z,W$ on $\tilde{M}$. A $3$-dimensional
non-Sasakian $(\kappa ,\mu )$-manifold has a constant $\varphi $-sectional
curvature, but for higher dimension this is not in general true. A
non-Sasakian $(\kappa ,\mu )$-manifold is of constant $\varphi $-sectional
curvature $c$ if and only if $\mu =\kappa +1$; in this case $c=-2k-1$. In
particular, the tangent sphere bundle of a manifold of constant curvature $%
c\neq 1$ has $\varphi $-sectional curvature $c^{2}=9\pm 4\sqrt{5}$ if and
only if $c=2\pm \sqrt{5}$. For more details we refer to \cite{Blair-book},
\cite{BKP-95} and \cite{Koufog-97}.

\section{\protect\boldmath$C$-totally real warped product submanifolds\label%
{sect-C-tot-real-wp-submfd}}

A submanifold $M$ in a contact manifold is called a $C${\em -totally real
submanifold} \cite{YKI-76} if every tangent vector of $M$ belongs to the
contact distribution. Thus, a submanifold $M$ in a contact metric manifold
is a $C${\em -totally real submanifold} if $\xi $ is normal to $M$. A
submanifold $M$ in an almost contact metric manifold is called {\em %
anti-invariant} \cite{YK-76} if $\varphi \left( TM\right) \subset T^{\perp
}M $. If a submanifold $M$ in a contact metric manifold is normal to the
structure vector field $\xi $, then it is anti-invariant. Thus $C$-totally
real submanifolds in a contact metric manifold are anti-invariant, as they
are normal to $\xi $. \medskip

\begin{ex-new}
\label{ex-C-tot-real-submfd} Given a point $p$ of a non-Sasakian $(\kappa
,\mu )$-manifold, there are at least two $C$-totally real submanifolds,
passing through $p$. To see this, consider the foliations given by the
eigendistributions of $h$; then their leaves are totally geodesic $C$%
-totally real submanifolds of the given non-Sasakian $(\kappa ,\mu )$%
-manifold.
\end{ex-new}

For a $C$-totally real submanifold in a contact metric manifold we have
\[
\left\langle A_{\xi }X,Y\right\rangle =\left\langle -\tilde{\nabla}_{X}\xi
,Y\right\rangle =\left\langle \varphi X+\varphi hX,Y\right\rangle ,
\]
which implies that
\begin{equation}
A_{\xi }=\left( \varphi h\right) ^{T},  \label{eq-A-xi}
\end{equation}
where $\left( \varphi h\right) ^{T}X$ is the tangential part of $\varphi hX$
for all $X\in TM$.

\medskip Now, we obtain a basic inequality involving the Laplacian of the
warping function and the squared mean curvature of a $C$-totally real warped
product submanifold of a $\left( \kappa ,\mu \right) $-space form.

\begin{th}
\label{th-laplacian-f-H-kmu-sp-form} Let $x$ be a $C$-totally real immersion
of an $n$-dimensional warped product manifold $M_{1}\times _{f}M_{2}$ into a
$\left( 2m+1\right) $-dimensional $\left( \kappa ,\mu \right) $-space form $%
\tilde{M}\left( c\right) $. Then
\begin{eqnarray}
\frac{\Delta f}{f} &\leq &\frac{n^{2}}{4n_{2}}\,\Vert H\Vert ^{2}+\frac{1}{4}%
\,n_{1}\left( c+3\right) +{\rm trace}(h^{T}|_{M_{1}})+\frac{n_{1}}{n_{2}}\,%
{\rm trace}(h^{T}|_{M_{2}})  \nonumber \\
&&+\,\frac{1}{4n_{2}}\left\{ ({\rm trace}(h^{T}))^{2}-({\rm trace}%
(h^{T}|_{M_{1}}))^{2}-({\rm trace}(h^{T}|_{M_{2}}))^{2}\right.   \nonumber \\
&&\left. -\,({\rm trace}(A_{\xi }))^{2}+({\rm trace}(A_{\xi
}|_{M_{1}}))^{2}+({\rm trace}(A_{\xi }|_{M_{2}}))^{2}\right.   \nonumber \\
&&\left. -\,\Vert h^{T}\Vert ^{2}+\Vert h^{T}|_{M_{1}}\Vert ^{2}+\Vert
h^{T}|_{M_{2}}\Vert ^{2}\right. \left. +\,\Vert A_{\xi }\Vert ^{2}-\Vert
A_{\xi }|_{M_{1}}\Vert ^{2}-\Vert A_{\xi }|_{M_{2}}\Vert ^{2}\right\} ,
\label{eq-lap-f-H-kmu-perp-1}
\end{eqnarray}
where $n_{i}=\dim M_{i}$, $i=1,2$, and $\Delta $ is the Laplacian operator
of $M_{1}$. The equality case of $(\ref{eq-lap-f-H-kmu-perp-1})$ holds
identically if and only if $x$ is a mixed totally geodesic immersion and $%
n_{1}H_{1}=n_{2}H_{2}$, where $H_{i}$, $i=1,2,$ are the partial mean
curvature vectors.
\end{th}

\noindent {\bf Proof.} Let $M_{1}\times _{f}M_{2}$ be a $C$-totally real
warped product submanifold into a $\left( \kappa ,\mu \right) $-space form $%
\tilde{M}\left( c\right) $ of constant $\varphi $-sectional curvature $c$.
We choose a local orthonormal frame $\{e_{1},\ldots ,e_{n},e_{n+1},\ldots
,e_{2m+1}\}$ such that $e_{1},\ldots ,e_{n_{1}}$ are tangent to $M_{1}$, $%
e_{n_{1}+1},\ldots ,e_{n}$ are tangent to $M_{2}$. Then from (\ref
{eq-Riem-curv}) and (\ref{eq-A-xi}) we have
\begin{eqnarray}
\tilde{K}\left( e_{i}\wedge e_{j}\right)  &=&\frac{c+3}{4}+\langle
h^{T}e_{i},e_{i}\rangle +\langle h^{T}e_{j},e_{j}\rangle   \nonumber \\
&&+\frac{1}{2}\,\left\{ \langle h^{T}e_{i},e_{i}\rangle \langle
h^{T}e_{j},e_{j}\rangle -\langle h^{T}e_{i},e_{j}\rangle ^{2}\right.
\nonumber \\
&&\qquad \quad \left. -\,\langle A_{\xi }e_{i},e_{i}\rangle \langle A_{\xi
}e_{j},e_{j}\rangle +\langle A_{\xi }e_{i},e_{j}\rangle ^{2}\right\} ,
\label{eq-K-tilde-ij-1}
\end{eqnarray}
where $h^{T}X$ is the tangential part of $hX$ for $X\in TM$.

For a $k$-plane section $\Pi _{k}$, from (\ref{eq-K-tilde-ij-1}) it follows
that
\begin{eqnarray}
\tilde{\tau}\left( \Pi _{k}\right)  &=&{\frac{1}{8}}\,k\left( k-1\right)
\left( c+3\right) +\left( k-1\right) {\rm trace}(h^{T}|_{\Pi _{k}})
\nonumber \\
&&+\,\frac{1}{4}\,\left\{ ({\rm trace}(h^{T}|_{\Pi _{k}}))^{2}-\Vert
h^{T}|_{\Pi _{k}}\Vert ^{2}\right. \left. -\,({\rm trace}(A_{\xi }|_{\Pi
_{k}}))^{2}+\Vert A_{\xi }|_{\Pi _{k}}\Vert ^{2}\right\} .
\label{eq-tau-Pi-k-1}
\end{eqnarray}
Consequently, we have
\begin{eqnarray}
\tilde{\tau}\left( T_{p}M\right)  &=&{\frac{1}{8}}\,n\left( n-1\right)
\left( c+3\right) +\left( n-1\right) {\rm trace}\left( h^{T}\right)
\nonumber \\
&&+\ \frac{1}{4}\left\{ ({\rm trace}(h^{T}))^{2}-\Vert h^{T}\Vert
^{2}\right. \left. -({\rm trace}(A_{\xi }))^{2}+\Vert A_{\xi }\Vert
^{2}\right\} ,  \label{eq-tau-tilde-1}
\end{eqnarray}
\begin{eqnarray}
\tilde{\tau}\left( T_{p}M_{i}\right)  &=&{\frac{1}{8}}\,n_{i}\left(
n_{i}-1\right) \left( c+3\right) +\left( n_{i}-1\right) {\rm trace}%
(h^{T}|_{M_{i}})  \nonumber \\
&&+\,\frac{1}{4}\,\left\{ ({\rm trace}(h^{T}|_{M_{i}}))^{2}-\Vert
h^{T}|_{M_{i}}\Vert ^{2}\right.   \nonumber \\
&&\left. -\,({\rm trace}(A_{\xi }|_{M_{i}}))^{2}+\Vert A_{\xi
}|_{M_{i}}\Vert ^{2}\right\} ,  \label{eq-tau-tilde-1i}
\end{eqnarray}
where $i=1,2$. Using (\ref{eq-tau-tilde-1}) and (\ref{eq-tau-tilde-1i}) in (%
\ref{eq-lap-f-H}) we get (\ref{eq-lap-f-H-kmu-perp-1}). $\blacksquare $
\medskip

Putting $h=0$ in (\ref{eq-lap-f-H-kmu-perp-1}), we have the following

\begin{cor}
{\rm (Theorem 3.1, \cite{Mat-Mih-02-SUT})} Let $x$ be a $\,C$-totally real
immersion of an $n$-dimensional warped product manifold $M_{1}\times
_{f}M_{2}$ into a Sasakian space form. Then
\begin{equation}
\frac{\Delta f}{f}\leq \frac{n^{2}}{4n_{2}}\Vert H\Vert ^{2}+\frac{1}{4}%
n_{1}\left( c+3\right) ,  \label{eq-lap-f-H-kmu-perp-Sas}
\end{equation}
where $n_{i}=\dim M_{i}$, $i=1,2$, and $\Delta $ is the Laplacian operator
of $M_{1}$. The equality case of $(\ref{eq-lap-f-H-kmu-perp-Sas})$ holds
identically if and only if $x$ is a mixed totally geodesic immersion and $%
n_{1}H_{1}=n_{2}H_{2}$, where $H_{i}$, $i=1,2,$ are the partial mean
curvature vectors.
\end{cor}

Next, we establish a sharp relationship between the warped function $f$ of a
$C$-totally real warped product submanifold $M_{1}\times _{f}M_{2}$
isometrically immersed in a non-Sasakian $\left( \kappa ,\mu \right) $%
-manifold $\tilde{M}$ and the squared mean curvature $\Vert H\Vert ^{2}$ in
the following

\begin{th}
Let $x$ be a $C$-totally real immersion of an $n$-dimensional warped product
manifold $M_{1}\times _{f}M_{2}$ into a $\left( 2m+1\right) $-dimensional
non-Sasakian $\left( \kappa ,\mu \right) $-manifold $\tilde{M}$. Then
\begin{eqnarray}
&&\frac{\Delta f}{f}\leq \frac{n^{2}}{4n_{2}}\Vert H\Vert ^{2}+n_{1}\left( 1-%
\frac{\mu }{2}\right) +{\rm trace}\left( h^{T}|_{M_{1}}\right) +\frac{n_{1}}{%
n_{2}}{\rm trace}\left( h^{T}|_{M_{2}}\right)   \nonumber \\
&&{+}\,\frac{1}{2n_{2}}\left( \frac{1-\left( \mu /2\right) }{1-\kappa }%
\right) \{({\rm trace}(h^{T}))^{2}-({\rm trace}(h^{T}|_{M_{1}}))^{2}-({\rm %
trace}(h^{T}|_{M_{2}}))^{2}\}  \nonumber \\
&&-\,\frac{1}{2n_{2}}\left( \frac{\kappa -\left( \mu /2\right) }{1-\kappa }%
\right) \{({\rm trace}(A_{\xi }))^{2}-({\rm trace}(A_{\xi }|_{M_{1}}))^{2}-(%
{\rm trace}(A_{\xi }|_{M_{2}}))^{2}\}  \nonumber \\
&&-\,\frac{1}{2n_{2}}\left( \frac{1-\left( \mu /2\right) }{1-\kappa }\right)
\left\{ \Vert h^{T}\Vert ^{2}-\Vert h^{T}|_{M_{1}}\Vert ^{2}-\Vert
h^{T}|_{M_{2}}\Vert ^{2}\right\}   \nonumber \\
&&+\,\frac{1}{2n_{2}}\left( \frac{\kappa -\left( \mu /2\right) }{1-\kappa }%
\right) \{\Vert A_{\xi }\Vert ^{2}-\Vert A_{\xi }|_{M_{1}}\Vert ^{2}-\Vert
A_{\xi }|_{M_{2}}\Vert ^{2}\},  \label{eq-lap-f-H-kmu-perp-2}
\end{eqnarray}
where $n_{i}=\dim M_{i}$, $i=1,2$, and $\Delta $ is the Laplacian operator
of $M_{1}$. The equality case of $(\ref{eq-lap-f-H-kmu-perp-2})$ holds
identically if and only if $x$ is a mixed totally geodesic immersion and $%
n_{1}H_{1}=n_{2}H_{2}$, where $H_{i}$, $i=1,2,$ are the partial mean
curvature vectors.
\end{th}

\noindent {\bf Proof.} Let $M_{1}\times _{f}M_{2}$ be a $C$-totally real
warped product submanifold of a non-Sasakian $(\kappa ,\mu )$-manifold $%
\tilde{M}$. We choose a local orthonormal frame $\{e_{1},\ldots
,e_{n},\allowbreak e_{n+1},\ldots ,e_{2m+1}\}$ such that $e_{1},\ldots
,e_{n_{1}}$ are tangent to $M_{1}$, $e_{n_{1}+1},\ldots ,e_{n}$ are tangent
to $M_{2}$. Then from (\ref{eq-kmu-curvature}) we have
\begin{eqnarray}
&&\tilde{K}\left( e_{i}\wedge e_{j}\right) =\left( 1-\frac{\mu }{2}\right)
+\left\langle h^{T}e_{i},e_{i}\right\rangle +\left\langle
h^{T}e_{j},e_{j}\right\rangle  \nonumber \\
&&\qquad +\,\frac{1-\left( \mu /2\right) }{1-\kappa }\left\{ \left\langle
h^{T}e_{i},e_{i}\right\rangle \left\langle h^{T}e_{j},e_{j}\right\rangle
-\left\langle h^{T}e_{i},e_{j}\right\rangle ^{2}\right\}  \nonumber \\
&&\qquad +\,\frac{\kappa -\left( \mu /2\right) }{1-\kappa }\left\{ \langle
A_{\xi }e_{i},e_{i}\rangle \langle A_{\xi }e_{j},e_{j}\rangle -\langle
A_{\xi }e_{i},e_{j}\rangle ^{2}\right\} .  \label{eq-K-tilde-ij-2}
\end{eqnarray}
For a $k$-plane section $\Pi _{k}$, from (\ref{eq-K-tilde-ij-2}) we obtain
\begin{eqnarray}
\tilde{\tau}\left( \Pi _{k}\right) &=&{\frac{1}{2}}\,k\left( k-1\right)
\left( 1-\frac{\mu }{2}\right) +\left( k-1\right) {\rm trace}(h^{T}|_{\Pi
_{k}})  \nonumber \\
&&{+}\,\frac{1}{2}\left( \frac{1-\left( \mu /2\right) }{1-\kappa }\right)
\left\{ ({\rm trace}(h^{T}|_{\Pi _{k}}))^{2}-\Vert h^{T}|_{\Pi _{k}}\Vert
^{2}\right\}  \nonumber \\
&&-\,\frac{1}{2}\left( \frac{\kappa -\left( \mu /2\right) }{1-\kappa }%
\right) \left\{ ({\rm trace}(A_{\xi }|_{\Pi _{k}}))^{2}-\Vert A_{\xi }|_{\Pi
_{k}}\Vert ^{2}\right\} .  \label{eq-tau-Pi-k-2}
\end{eqnarray}
Consequently, we have
\begin{eqnarray}
\tilde{\tau}\left( T_{p}M\right) &=&{\frac{1}{2}}\,n\left( n-1\right) \left(
1-\frac{\mu }{2}\right) +\left( n-1\right) {\rm trace}(h^{T})  \nonumber \\
&&{+}\,\frac{1}{2}\,\left( \frac{1-\left( \mu /2\right) }{1-\kappa }\right)
\left\{ ({\rm trace}(h^{T}))^{2}-\Vert h^{T}\Vert ^{2}\right\}  \nonumber \\
&&-\,\frac{1}{2}\,\left( \frac{\kappa -\left( \mu /2\right) }{1-\kappa }%
\right) \left\{ ({\rm trace}(A_{\xi }))^{2}-\Vert A_{\xi }\Vert ^{2}\right\}
,  \label{eq-tau-tilde-2}
\end{eqnarray}
\begin{eqnarray}
\tilde{\tau}\left( T_{p}M_{i}\right) &=&{\frac{1}{2}}\,n_{i}\left(
n_{i}-1\right) \left( 1-\frac{\mu }{2}\right) +\left( n_{i}-1\right) {\rm %
trace}(h^{T}|_{M_{i}})  \nonumber \\
&&{+}\,\frac{1}{2}\,\left( \frac{1-\left( \mu /2\right) }{1-\kappa }\right)
\left\{ ({\rm trace}(h^{T}|_{M_{i}}))^{2}-\Vert h^{T}|_{M_{i}}\Vert
^{2}\right\}  \nonumber \\
&&-\,\frac{1}{2}\,\left( \frac{\kappa -\left( \mu /2\right) }{1-\kappa }%
\right) \left\{ ({\rm trace}(A_{\xi }|_{M_{i}}))^{2}-\Vert A_{\xi
}|_{M_{i}}\Vert ^{2}\right\} ,  \label{eq-tau-tilde-2i}
\end{eqnarray}
where $i=1,2$. Putting values from (\ref{eq-tau-tilde-2}) and (\ref
{eq-tau-tilde-2i}) in (\ref{eq-lap-f-H}) we obtain (\ref
{eq-lap-f-H-kmu-perp-2}). $\blacksquare $ \medskip

As applications, we derive certain obstructions to the existence of minimal $%
C$-totally real warped product submanifolds in $\left( \kappa ,\mu \right) $%
-space forms, non-Sasakian $\left( \kappa ,\mu \right) $-manifolds and
Sasakian space forms. \medskip

\begin{cor}
\label{cor-appl-1} Let $M_{1}\times _{f}M_{2}$ be a warped product manifold,
whose warping function $f$ is harmonic. Then:

\begin{enumerate}
\item[{\bf (a)}]  $M_{1}\times _{f}M_{2}$ can not admit minimal $C$-totally
real immersion into a $\left( \kappa ,\mu \right) $-space form $\tilde{M}%
\left( c\right) $ with
\begin{eqnarray}
0 &>&\frac{1}{4}\,n_{1}\left( c+3\right) +{\rm trace}\left(
h^{T}|_{M_{1}}\right) +\frac{n_{1}}{n_{2}}\,{\rm trace}\left(
h^{T}|_{M_{2}}\right)   \nonumber \\
&&+\,\frac{1}{4n_{2}}\left\{ \left( {\rm trace}\left( h^{T}\right) \right)
^{2}-\left( {\rm trace}\left( h^{T}|_{M_{1}}\right) \right) ^{2}-\left( {\rm %
trace}\left( h^{T}|_{M_{2}}\right) \right) ^{2}\right.   \nonumber \\
&&\left. -\,({\rm trace}(A_{\xi }))^{2}+({\rm trace}(A_{\xi
}|_{M_{1}}))^{2}+({\rm trace}(A_{\xi }|_{M_{2}})^{2}\right.   \nonumber \\
&&\left. -\,\left\| h^{T}\right\| ^{2}+\left\| h^{T}|_{M_{1}}\right\|
^{2}+\left\| h^{T}|_{M_{2}}\right\| ^{2}\right. \left. +\,\left\| A_{\xi
}\right\| ^{2}-\left\| A_{\xi }|_{M_{1}}\right\| ^{2}-\left\| A_{\xi
}|_{M_{2}}\right\| ^{2}\right\} .  \label{eq-cor-appl-1-1}
\end{eqnarray}

\item[{\bf (b)}]  Every minimal $C$-totally real immersion of $M_{1}\times
_{f}M_{2}$ in a $\left( \kappa ,\mu \right) $-space form $\tilde{M}\left(
c\right) $ with
\begin{eqnarray}
0 &=&\frac{1}{4}\,n_{1}\left( c+3\right) +{\rm trace}\left(
h^{T}|_{M_{1}}\right) +\frac{n_{1}}{n_{2}}\,{\rm trace}\left(
h^{T}|_{M_{2}}\right)   \nonumber \\
&&+\,\frac{1}{4n_{2}}\left\{ \left( {\rm trace}\left( h^{T}\right) \right)
^{2}-\left( {\rm trace}\left( h^{T}|_{M_{1}}\right) \right) ^{2}-\left( {\rm %
trace}\left( h^{T}|_{M_{2}}\right) \right) ^{2}\right.   \nonumber \\
&&\left. -\,({\rm trace}(A_{\xi }))^{2}+({\rm trace}(A_{\xi
}|_{M_{1}}))^{2}+({\rm trace}(A_{\xi }|_{M_{2}}))^{2}\right.   \nonumber \\
&&\left. -\,\left\| h^{T}\right\| ^{2}+\left\| h^{T}|_{M_{1}}\right\|
^{2}+\left\| h^{T}|_{M_{2}}\right\| ^{2}\right. \left. +\,\left\| A_{\xi
}\right\| ^{2}-\left\| A_{\xi }|_{M_{1}}\right\| ^{2}-\left\| A_{\xi
}|_{M_{2}}\right\| ^{2}\right\}   \label{eq-cor-appl-1-2}
\end{eqnarray}
is a warped product immersion.
\end{enumerate}
\end{cor}

\noindent {\bf Proof.} We assume that the warping function $f$ is a harmonic
function on $M_{1}$ and $M_{1}\times _{f}M_{2}$ admits a minimal $C$-totally
real immersion in a $\left( \kappa ,\mu \right) $-space form $\tilde{M}%
\left( c\right) $. Then, the inequality (\ref{eq-lap-f-H-kmu-perp-1})
becomes
\begin{eqnarray*}
0 &\leq &\frac{1}{4}n_{1}\left( c+3\right) +{\rm trace}\left(
h^{T}|_{M_{1}}\right) +\frac{n_{1}}{n_{2}}{\rm trace}\left(
h^{T}|_{M_{2}}\right) \\
&&+\,\frac{1}{4n_{2}}\left\{ \left( {\rm trace}\left( h^{T}\right) \right)
^{2}-\left( {\rm trace}\left( h^{T}|_{M_{1}}\right) \right) ^{2}-\left( {\rm %
trace}\left( h^{T}|_{M_{2}}\right) \right) ^{2}\right. \\
&&\left. -\,({\rm trace}(A_{\xi }))^{2}+({\rm trace}(A_{\xi
}|_{M_{1}}))^{2}+({\rm trace}(A_{\xi }|_{M_{2}}))^{2}\right. \\
&&\left. -\,\left\| h^{T}\right\| ^{2}+\left\| h^{T}|_{M_{1}}\right\|
^{2}+\left\| h^{T}|_{M_{2}}\right\| ^{2}\right. \left. +\left\| A_{\xi
}\right\| ^{2}-\left\| A_{\xi }|_{M_{1}}\right\| ^{2}-\left\| A_{\xi
}|_{M_{2}}\right\| ^{2}\right\} .
\end{eqnarray*}
This proves {\bf (a)}.

Now, we prove {\bf (b)}. If (\ref{eq-cor-appl-1-2}) is true then the
equality case of (\ref{eq-lap-f-H-kmu-perp-1}) is true and by Theorem~\ref
{th-laplacian-f-H-kmu-sp-form}, it follows that $M_{1}\times _{f}M_{2}$ is
mixed totally geodesic and $H_{1}=H_{2}=0$. Then a well-known result of
N\"{o}lker \cite{Nolker-96} implies that the immersion is a warped product
immersion. $\blacksquare $ \medskip

Similar to the above Corollary, we have

\begin{cor}
\label{cor-appl-2} Let $M_{1}\times _{f}M_{2}$ be a warped product manifold,
whose warping function $f$ is harmonic. Then:

\begin{enumerate}
\item[{\bf (a)}]  $M_{1}\times _{f}M_{2}$ can not admit minimal $C$-totally
real immersion into a non-Sasakian $\left( \kappa ,\mu \right) $-manifold
with
\begin{eqnarray}
0 &>&n_{1}\left( 1-\frac{\mu }{2}\right) +{\rm trace}\left(
h^{T}|_{M_{1}}\right) +\frac{n_{1}}{n_{2}}\,{\rm trace}\left(
h^{T}|_{M_{2}}\right)   \nonumber \\
&&{+}\,\frac{1}{2n_{2}}\left( \frac{1-\left( \mu /2\right) }{1-\kappa }%
\right) \{({\rm trace}(h^{T}))^{2}-({\rm trace}(h^{T}|_{M_{1}}))^{2}-({\rm %
trace}(h^{T}|_{M_{2}}))^{2}\}  \nonumber \\
&&-\,\frac{1}{2n_{2}}\left( \frac{\kappa -\left( \mu /2\right) }{1-\kappa }%
\right) \{({\rm trace}(A_{\xi }))^{2}-({\rm trace}(A_{\xi }|_{M_{1}}))^{2}-(%
{\rm trace}(A_{\xi }|_{M_{2}}))^{2}\}  \nonumber \\
&&-\,\frac{1}{2n_{2}}\left( \frac{1-\left( \mu /2\right) }{1-\kappa }\right)
\left\{ \Vert h^{T}\Vert ^{2}-\Vert h^{T}|_{M_{1}}\Vert ^{2}-\Vert
h^{T}|_{M_{2}}\Vert ^{2}\right\}   \nonumber \\
&&+\,\frac{1}{2n_{2}}\left( \frac{\kappa -\left( \mu /2\right) }{1-\kappa }%
\right) \{\Vert A_{\xi }\Vert ^{2}-\Vert A_{\xi }|_{M_{1}}\Vert ^{2}-\Vert
A_{\xi }|_{M_{2}}\Vert ^{2}\}.  \label{eq-cor-appl-2-1}
\end{eqnarray}

\item[{\bf (b)}]  Every minimal $C$-totally real immersion of $M_{1}\times
_{f}M_{2}$ in a non-Sasakian $\left( \kappa ,\mu \right) $-manifold with
\begin{eqnarray}
0 &=&n_{1}\left( 1-\frac{\mu }{2}\right) +{\rm trace}\left(
h^{T}|_{M_{1}}\right) +\frac{n_{1}}{n_{2}}\,{\rm trace}\left(
h^{T}|_{M_{2}}\right)   \nonumber \\
&&{+}\,\frac{1}{2n_{2}}\left( \frac{1-\left( \mu /2\right) }{1-\kappa }%
\right) \{({\rm trace}(h^{T}))^{2}-({\rm trace}(h^{T}|_{M_{1}}))^{2}-({\rm %
trace}(h^{T}|_{M_{2}}))^{2}\}  \nonumber \\
&&-\,\frac{1}{2n_{2}}\left( \frac{\kappa -\left( \mu /2\right) }{1-\kappa }%
\right) \{({\rm trace}(A_{\xi }))^{2}-({\rm trace}(A_{\xi }|_{M_{1}}))^{2}-(%
{\rm trace}(A_{\xi }|_{M_{2}}))^{2}\}  \nonumber \\
&&-\,\frac{1}{2n_{2}}\left( \frac{1-\left( \mu /2\right) }{1-\kappa }\right)
\left\{ \Vert h^{T}\Vert ^{2}-\Vert h^{T}|_{M_{1}}\Vert ^{2}-\Vert
h^{T}|_{M_{2}}\Vert ^{2}\right\}   \nonumber \\
&&+\,\frac{1}{2n_{2}}\left( \frac{\kappa -\left( \mu /2\right) }{1-\kappa }%
\right) \{\Vert A_{\xi }\Vert ^{2}-\Vert A_{\xi }|_{M_{1}}\Vert ^{2}-\Vert
A_{\xi }|_{M_{2}}\Vert ^{2}\}  \label{eq-cor-appl-2-2}
\end{eqnarray}
is a warped product immersion.
\end{enumerate}
\end{cor}

Using $h=0$ in Corollary~\ref{cor-appl-1}, we immediately get the following

\begin{cor}
\label{cor-appl-1-Sas} {\rm \cite[Corollary 3.2]{Mat-Mih-02-SUT}} Let $%
M_{1}\times _{f}M_{2}$ be a warped product manifold, whose warping function $%
f$ is harmonic. Then:

\begin{enumerate}
\item[{\bf (a)}]  $M_{1}\times _{f}M_{2}$ can not admit minimal $C$-totally
real immersion into a Sasakian space form $\tilde{M}(c)$ with $c<-3$.

\item[{\bf (b)}]  Every minimal $C$-totally real immersion of $M_{1}\times
_{f}M_{2}$ in the standard Sasakian space form ${\Bbb R}^{2m+1}(-3)$ is a
warped product immersion.
\end{enumerate}
\end{cor}

We also have the following

\begin{cor}
{\rm \cite[Corollary 3.3]{Mat-Mih-02-SUT}} If the warping function $f$ of a
warped product $M_{1}\times _{f}M_{2}$ is an eigenfunction of the Laplacian
on $M_{1}$ with corresponding eigenvalue $\lambda >0$, then $M_{1}\times
_{f}M_{2}$ does not admit a minimal $C$-totally real immersion in a Sasakian
space form $\tilde{M}(c)$ with $c\leq -3$.
\end{cor}

In the following we obtain an example of a warped product $C$-totally real
submanifold of a non-Sasakian $(\kappa ,\mu )$-manifold which satisfies the
equality case of (\ref{eq-lap-f-H-kmu-perp-2}).

\begin{ex-new}
Let $(\tilde{M},\tilde{g})$ be a $4$-dimensional Riemannian manifold of
constant sectional curvature $c\neq 1$. Then its tangent sphere bundle $T_{1}%
\tilde{M}$ with the standard contact metric structure is a non-Sasakian $%
(\kappa ,\mu )$-manifold with $\kappa =c\left( 2-c\right) ,\ \mu =2c$ \cite
{BKP-95}. Let $M$ be a totally geodesic hypersurface of $\tilde{M}$. Then $M$
equipped with the induced Riemannian metric $g$ has constant sectional
curvature $c$ and also its tangent sphere bundle $T_{1}M$ is a $(\kappa ,\mu
)$-manifold with $\kappa =c(2-c),\ \mu =2c$. From Example~3.3 of \cite
{CTT-08} it follows that $T_{1}M$ is an invariant submanifold of $T_{1}%
\tilde{M}$. Let $S$ be a minimal $C$-totally real surface of $T_{1}M$ (which
is always possible in view of Example~\ref{ex-C-tot-real-submfd}). Define
the warped product manifold
\[
\left( -\frac{\pi }{2}\,,\frac{\pi }{2}\right) \times _{\cos t}S.
\]
Then the immersion $x:\left( -\frac{\pi }{2}\,,\frac{\pi }{2}\right) \times
_{\cos t}S\;\rightarrow \;T_{1}\tilde{M}$ defined by
\[
x\left( t,p\right) =\left( \sin t\right) N+\left( \cos t\right) p\,,
\]
where $N$ is a unit vector orthogonal to the linear subspace containing $%
T_{1}M$, is a $C$-totally real isometric immersion and satisfies the
equality case of (\ref{eq-lap-f-H-kmu-perp-2}).
\end{ex-new}

\noindent Department of Mathematics

\noindent Banaras Hindu University

\noindent Varanasi, 221 005, India

\noindent Email: {\tt mmtripathi66@yahoo.com}

\end{document}